\documentclass[10pt,letterpaper,conference]{ieeeconf}

\usepackage{latexsym, amssymb, amsmath}
\usepackage[dvips]{graphicx} 
\usepackage{enumerate}
\usepackage{flushend} 
\usepackage{mathrsfs} 
\usepackage{stfloats}
\usepackage{cite}

\allowdisplaybreaks[3] 

\newtheorem{prop}{Proposition} 
\newtheorem{proposition}[prop]{Proposition} 
\newtheorem{lem}{Lemma} 
\newtheorem{lemma}[lem]{Lemma} 
\newtheorem{thm}{Theorem} 
\newtheorem{theorem}[thm]{Theorem} 
\newtheorem{cor}{Corollary} 
\newtheorem{corollary}[cor]{Corollary} 
\newtheorem{defn}{Definition} 
\newtheorem{definition}[defn]{Definition} 
\newtheorem{exmp}{Example} 
\newtheorem{example}[exmp]{Example} 
\newtheorem{exam*}{Example}

\def\custombibliography#1{
 \normalsize
\section*{\centering References}
 \list
 {[\arabic{enumi}]}{\settowidth\labelwidth{[#1]}\leftmargin\labelwidth
 \setlength{\itemsep}{.1em}
 \advance\leftmargin\labelsep
 \usecounter{enumi}}
 \def\newblock{\hskip .11em plus .33em minus -.07em}
 \sloppy
 \sfcode`\.=1000\relax}

\def\L2{{\cal L}_2}
\def\begar{\begin{array}}
\def\endar{\end{array}}
\def\begce{\begin{center}}
\def\endce{\end{center}}
\def\begco{\begin{cor}}
\def\endco{\end{cor}}
\def\begde{\begin{defn}}
\def\endde{\end{defn}}
\def\begdes{\begin{description}}
\def\enddes{\end{description}}
\def\begdi{\begin{displaymath}}
\def\enddi{\end{displaymath}}
\def\begdis{\begin{eqnarray*}}
\def\enddis{\end{eqnarray*}}
\def\begen{\begin{enumerate}}
\def\enden{\end{enumerate}}
\def\begeq{\begin{equation}}
\def\endeq{\end{equation}}
\def\begeqa{\begin{eqnarray}}
\def\endeqa{\end{eqnarray}}
\def\begex{\begin{exmp}}
\def\endex{\end{exmp}}
\def\begfig{\begin{fig}}
\def\endfig{\end{fig}}
\def\begit{\begin{itemize}}
\def\endit{\end{itemize}}
\def\begle{\begin{lem}}
\def\endle{\end{lem}}
\def\begpro{\begin{prop}} 
\def\endpro{\end{prop}} 
\def\begth{\begin{thm}}
\def\endth{\end{thm}}

\def\begres{\noindent{\bf Remarks}:\begin{enumerate}}
\def\endres{\end{enumerate} \par}
\def\begpr{\begin{proof}} 
\def\endpr{\end{proof}}

\newcommand\bull{\vrule height .9ex width .8ex depth -.1ex } 

\newcommand\re{\rm I\! R}
\newcommand\cdcout[1]{} 

\newcommand{\rv}[1]{\boldsymbol{#1}} 
\newcommand{\RomanNumber}[1]{\uppercase\expandafter{\romannumeral #1}}
\newcommand{\romannumber}[1]{\lowercase\expandafter{\romannumeral #1}}

\DeclareMathAlphabet{\mathpzc}{OT1}{pzc}{m}{it}

\def\1{\rv 1} 

\renewcommand{\baselinestretch}{0.94}
\IEEEoverridecommandlockouts                            
\title{Collision Avoidance for Bi-Steerable Car Using Analytic Left Inversion}

\author{Luis A. Duffaut Espinosa$^\dag$\thanks{$^\dag$Department of Electrical
and Computer Engineering, George Mason University, Fairfax, VA 22030, USA.}
$\;\;\;$
W.~Steven Gray$^\ddag$\thanks{$^\ddag$Department
of Electrical and Computer Engineering, Old Dominion University, Norfolk,
Virginia 23529-0246, USA.}
}

\usepackage{cite} 
\usepackage{color} 
\usepackage{graphicx} 
\usepackage{latexsym,amssymb,amsmath} 
\usepackage{mathrsfs} 
\usepackage{multicol}
\usepackage{fixltx2e} 

\allowdisplaybreaks[4]
\centerfigcaptionstrue

\definecolor{Light}{gray}{0.85}

\def\abs#1{\left\vert #1 \right\vert}
\def\allpolyx0degn{\mbox{$P_n$}}

\def\allseries{\mbox{$\re\langle\langle X \rangle\rangle$}}
\def\allseriesdelta{\mbox{$\re\langle\langle X_\delta \rangle\rangle$}}
\def\allseriesell{\mbox{$\re^{\ell} \langle\langle X \rangle\rangle$}}

\def\allseriesm{\mbox{$\re^m\langle\langle X \rangle\rangle$}}
\def\allseriesmtimesm{\mbox{$\re^{m\times m}\langle\langle X \rangle\rangle$}}

\def\allseriesellLC{\mbox{$\re^{\ell}_{LC}\langle\langle X \rangle\rangle$}}
\def\allseriesellGC{\mbox{$\re^{\ell}_{GC}\langle\langle X \rangle\rangle$}}

\def\allseriesXO{\mbox{$\re [[ X_0 ]]$}}

\def\allseriesXOmLC{\mbox{$\re^m_{LC} [[ X_0 ]]$}}

\def\allseriesX1{\mbox{$\re [[ X_1 ]]$}}
\def\bfem#1{{\bf \em #1}} 

\def\bull{\rule{0.08in}{0.08in}} 

\newcommand{\comment}[1]{} 

\def\eqref#1{(\ref{#1})} 

\def\Fliessdelta{\mathscr{F}_{\delta}}

\def\modcomp{\:\tilde{\circ}\,} 

\def\norm#1{\Vert#1\Vert}

\def\openbull{\framebox[0.08in][c]{$\;$}} 

\def\re{{\mathbb R}} 

\def\sameau{\rule[0.017in]{0.2in}{0.012in}}

\def\shuffle{{\scriptscriptstyle \;\sqcup \hspace*{-0.05cm}\sqcup\;}}

\def\supp{{\rm supp}}

\def\begals{\[\begin{aligned}}
\def\endals{\end{aligned}\]}
\def\begce{\begin{center}}
\def\endce{\end{center}}
\def\begar{\begin{array}}
\def\endar{\end{array}}
\def\begeq{\begin{equation}}
\def\endeq{\end{equation}}
\def\begdi{\begin{displaymath}}
\def\enddi{\end{displaymath}}
\def\begdis{\begin{eqnarray*}}
\def\enddis{\end{eqnarray*}}
\def\begeqa{\begin{eqnarray}}
\def\endeqa{\end{eqnarray}}
\def\begdes{\begin{description}}
\def\enddes{\end{description}}
\def\begit{\begin{itemize}}
\def\endit{\end{itemize}}
\def\begen{\begin{enumerate}}
\def\enden{\end{enumerate}}
\def\beglar{\left[\begin{array}}
\def\endrar{\end{array}\right]}
\def\begle{\begin{mylemma}}
\def\endle{\end{mylemma}}
\def\begde{\begin{mydefinition}}
\def\endde{\end{mydefinition}}
\def\begth{\begin{mytheorem}}
\def\endth{\end{mytheorem}}
\def\begco{\begin{mycorollary}}
\def\endco{\end{mycorollary}}
\def\begprop{\begin{myproposition}}
\def\endprop{\end{myproposition}}
\def\begex{\begin{myexample}}
\def\endex{\hfill\openbull \end{myexample} \vspace*{0.15in}}
\def\begexer{\begin{myexercise}}
\def\endexer{\end{myexercise}}

\def\begres{\noindent{\bf Remarks}:\begin{enumerate}}
\def\endres{\end{enumerate} \par}
\def\begpr{\noindent{\em Proof:}$\;\;$}
\def\endpr{\hfill\bull \vspace*{0.15in}}
\def\begtab{\begin{tabular}}
\def\endtab{\end{tabular}}
\def\rref#1{(\ref{#1})}
\def\allseriesA{\mbox{$\re\!\ll\!A\!\gg$}}
\def\allseriesA'{\mbox{$\re\!\ll\!A'\!\gg$}}

\def\co{{\rm co}}
\def\si{{\rm si}}

\def\begce{\begin{center}}
\def\endce{\end{center}}
\def\begar{\begin{array}}
\def\endar{\end{array}}
\def\begeq{\begin{equation}}
\def\endeq{\end{equation}}
\def\begdi{\begin{displaymath}}
\def\enddi{\end{displaymath}}
\def\begdis{\begin{eqnarray*}}
\def\enddis{\end{eqnarray*}}
\def\begeqa{\begin{eqnarray}}
\def\endeqa{\end{eqnarray}}
\def\begdes{\begin{description}}
\def\enddes{\end{description}}
\def\begit{\begin{itemize}}
\def\endit{\end{itemize}}
\def\begen{\begin{enumerate}}
\def\enden{\end{enumerate}}
\def\beglar{\left[\begin{array}}
\def\endrar{\end{array}\right]}
\def\begle{\begin{lemma}}
\def\endle{\end{lemma}}
\def\begde{\begin{definition}}
\def\endde{\end{definition}}
\def\begth{\begin{theorem}}
\def\endth{\end{theorem}}
\def\begco{\begin{corollary}}
\def\endco{\end{corollary}}
\def\begprop{\begin{proposition}}
\def\endprop{\end{proposition}}
\def\begex{\begin{example}}
\def\endex{\hfill\openbull \end{example} \vspace*{0.1in}}
\def\begexer{\begin{exercise}}
\def\endexer{\end{exercise}}

\def\begres{\noindent{\bf Remarks}:\begin{enumerate}}
\def\endres{\end{enumerate} \par}
\def\begpr{\noindent{\em Proof:}$\;\;$}
\def\endpr{\hfill\bull \vspace*{0.1in}}
\def\begtab{\begin{tabular}}
\def\endtab{\end{tabular}}
\def\rref#1{(\ref{#1})}

\date{\today}

\begin{document}

\maketitle

\begin{abstract}

A case study is presented of a collision avoidance
system that directly integrates the kinematics of a bi-steerable car with
a suitable path planning algorithm. The first step is to
identify a path using the method of rapidly exploring random trees, and then
a spline approximation is computed.
The second step is to solve the output tracking problem
by explicitly computing the left inverse of the kinematics of the system
to render the Taylor series of the desired input for each polynomial section
of the spline approximation. The method is demonstrated by numerical
simulation.

\end{abstract}

\vspace*{0.08in}
\begin{keywords}
Collision avoidance, motion planning, autonomous vehicles, system inversion, Chen-Fliess series
\end{keywords}

\section{Introduction}

A convergence of technological advances in conjunction with changing
societal attitudes suggests that autonomous wheeled vehicles are poised to soon
enter main stream applications in business, the government, and the private
realm \cite{Hars_10}. As this happens, it is likely that vehicle design will evolved in
new directions once the human driver and passengers are completely
replaced by a computer control system and cargo. For example, a
bi-steerable car with independent front and rear axels provides
improved maneuverability and handling, which would be invaluable in
urban settings to improve collision avoidance
\cite{Hermosillo-etal_03,Hermosillo-Sekhavat_03,Petrov_09,Sekhavat-et-al_2001,Sekhavat-etal_01a,Sekhavat-Hermosillo_02,Tchenderli-Braham-Hamerlain_13}.
But such improved performance is predicated on control methodologies that
properly combine the true dynamics of the vehicle with the
proper path planning tools
\cite{DeLuca-et-al_2001,Lavalle_2006}.

In this paper, a case study is presented of a collision avoidance
system that directly integrates the kinematics of a bi-steerable car with
a suitable path planning algorithm. Here collision avoidance means simply
steering a vehicle modeled as a point from a starting location to a final
location while avoiding all stationary obstacles in between the two locations.
The first step is solving the path planning problem. Numerous algorithms appear in
the literature which are applicable to the collision avoidance problem, such as
rapidly exploring random trees (RRT), probabilistic roadmap (PRM), artificial
potential fields, and genetic algorithms \cite{Lavalle_2006}.
The second step is solving the output tracking problem. That is, determine
the control inputs so that the vehicle accurately follows the desired path.
Mathematically, this corresponds to
computing a left inverse for the system. There is an extensive literature on
this topic in the context of geometric nonlinear control theory, see,
for example, \cite{Engl-et-al_05,Getz_95,Hirschorn_79,Hirschorn_79a,Isidori_95,Singh_81,Singh_82}.
This class of techniques generates left inverses dynamically by driving a
certain inverse dynamical system with the desired output. While very convenient in some
situations, this method has certain limitations, like requiring the system to be minimum phase. In
addition, this approach
does not provide an {\em explicit} representation of the control input one is seeking.
Another common approach is to use the property of {\em flatness}, which allows one
to explicitly compute a left inverse analytically provided suitable
flat outputs can be identified \cite{Fliess-et-al_95,Levine_2009}.
Unfortunately, these output are sometimes not the most desirable
from a physical point of view. This turns out to the case for the bi-steerable
car, where the flat outputs do not correspond to simple physical variables like
a center point on the vehicle or its speed or orientation
\cite{Hermosillo-Sekhavat_03,Sekhavat-et-al_2001,Sekhavat-etal_01a,Sekhavat-Hermosillo_02}.
This then makes integrating the tracking problem with path planning problem more complex.
A final inversion method, which is the one employed in this paper, is to use an analytic
expression for the left inverse which is available for systems which have a
well defined vector relative degree for the outputs of interest
\cite{Gray-etal_SCL14,Gray-etal_CDC15}. They do {\em not} have to be flat outputs. The method is based
on a Fliess operator representation of input-output system, and a
combinatorial Hopf algebra technique that renders an explicit formula for the Taylor series
of the left inverse when the output is an analytic function in the range of the
input-output map. This single formula can be pre-computed efficiently off-line to
arbitrary precision \cite{Duffaut-etal_JACxx,Gray-etal_MTNS14}, and,
for example, could be hardwired on an FPGA for real-time implementation. In this setting,
the desired output trajectory is approximated by a spline, and then
the left inverse of each polynomial section is numerically evaluated using
the inversion formula.
It will be shown here that this strategy provides a feasible and accurate solution to
the collision avoidance problem for a bi-steerable car. For brevity, only the
vehicle's kinematics are considered.

The paper is organized as follows. In Section \ref{sec:2}, preliminaries concerning
Fliess operators and their inverses are briefly summarized to make the presentation
more self-contained. Then, in Section
\ref{sec:left-inversion-Bisteered-system}, the left inverse of the kinematics
of the bi-steerable car is computed. This result
is then integrated in Section \ref{sec:simulations} with the RRT path planning algorithm
and the corresponding numerical
simulations are presented. Section \ref{sec:comclusions} provides the paper's conclusions.

\section{Preliminaries} \label{sec:2}

In this section, some preliminaries concerning the theory of Fliess operators are outlined
as they provide the cornerstone for the method presented in this paper.
The interested reader is referred to \cite{Gray-etal_Automatic_14,Gray-etal_SCL14}
for a more complete treatment.

\subsection{Fliess Operators and Their Interconnections}

A finite nonempty set of noncommuting symbols $X=\{ x_0,x_1,$ $\ldots,x_m\}$ is
called an {\em alphabet}. Each element of $X$ is called a {\em letter}, and any
finite sequence of letters from $X$, $\eta=x_{i_1}\cdots x_{i_k}$, is called a
{\em word} over $X$. The {\em length} of $\eta$, $\abs{\eta}$, is the number of
letters in $\eta$. The set of all words with length $k$ is denoted by
$X^k$. The set of all words including the empty word, $\emptyset$,
is designated by $X^\ast$. It forms a monoid under catenation.
The set $\eta X^\ast$ is comprised of all words with the prefix $\eta$.
Any mapping $c:X^\ast\rightarrow
\re^\ell$ is called a {\em formal power series}. The value of $c$ at
$\eta\in X^\ast$ is written as $(c,\eta)$ and called the {\em coefficient} of
$\eta$ in $c$.
Typically, $c$ is
represented as the formal sum $c=\sum_{\eta\in X^\ast}(c,\eta)\eta.$
If the {\em constant term} $(c,\emptyset)=0$ then $c$ is said to be {\em
proper}. The {\em support} of $c$, $\supp(c)$, is the set of all words having
nonzero coefficients. The collection of all formal power series over $X$ is
denoted by $\allseriesell$. It forms an associative $\re$-algebra under
the catenation product and a commutative and associative $\re$-algebra under
the shuffle product, denoted here by $\shuffle$. The latter is the
$\re$-bilinear extension of the shuffle product of two words, which is defined
inductively by
$(x_i\eta)\shuffle
(x_j\xi)=x_i(\eta\shuffle(x_j\xi))+x_j((x_i\eta)\shuffle \xi)$
with $\eta\shuffle\emptyset=\emptyset\shuffle\eta=\eta$ for all $\eta,\xi\in
X^\ast$ and $x_i,x_j\in X$.

One can formally associate with any series $c\in\allseriesell$ a causal
$m$-input, $\ell$-output operator, $F_c$, in the following manner.
Let $\mathfrak{p}\ge 1$ and $t_0 < t_1$ be given. For a Lebesgue measurable
function $u: [t_0,t_1] \rightarrow\re^m$, define
$\norm{u}_{\mathfrak{p}}=\max\{\norm{u_i}_{\mathfrak{p}}: \ 1\le
i\le m\}$, where $\norm{u_i}_{\mathfrak{p}}$ is the usual
$L_{\mathfrak{p}}$-norm for a measurable real-valued function,
$u_i$, defined on $[t_0,t_1]$.  Let $L^m_{\mathfrak{p}}[t_0,t_1]$
denote the set of all measurable functions defined on $[t_0,t_1]$
having a finite $\norm{\cdot}_{\mathfrak{p}}$ norm and
$B_{\mathfrak{p}}^m(R)[t_0,t_1]:=\{u\in
L_{\mathfrak{p}}^m[t_0,t_1]:\norm{u}_{\mathfrak{p}}\leq R\}$. Assume
$C[t_0,t_1]$ is the subset of continuous functions in $L_{1}^m[t_0,t_1]$.
Define inductively for each $\eta\in X^{\ast}$ the map $E_\eta:
L_1^m[t_0, t_1]\rightarrow C[t_0, t_1]$ by setting $E_\emptyset[u]=1$ and
letting \[E_{x_i\bar{\eta}}[u](t,t_0) =
\int_{t_0}^tu_{i}(\tau)E_{\bar{\eta}}[u](\tau,t_0)\,d\tau, \] where
$x_i\in X$, $\bar{\eta}\in X^{\ast}$, and $u_0=1$. The
input-output system corresponding to $c$ is the {\em Fliess operator}
\begdi
F_c[u](t) =
\sum_{\eta\in X^{\ast}} (c,\eta)\,E_\eta[u](t,t_0)
\label{eq:Fliess-operator-defined}.
\enddi
If there exist real numbers $K_c,M_c>0$ such that
$\abs{(c,\eta)}\le K_c M_c^{|\eta|}|\eta|!$, $\forall\eta\in X^{\ast}$,
then $F_c$ constitutes a well defined mapping from
$B_{\mathfrak p}^m(R)[t_0,$ $t_0+T]$ into $B_{\mathfrak
q}^{\ell}(S)[t_0, \, t_0+T]$ for sufficiently small $R,T >0$,
where the numbers $\mathfrak{p},\mathfrak{q}\in[1,\infty]$ are
conjugate exponents, i.e., $1/\mathfrak{p}+1/\mathfrak{q}=1$.
(Here, $\abs{z}:=\max_i \abs{z_i}$ when $z\in\re^\ell$.) The set of all such
{\em locally convergent} series is denoted by $\allseriesellLC$. On the other
hand, if $\abs{(c,\eta)}\le K_c M_c^{|\eta|}$, $\forall\eta\in X^{\ast}$, then the operator
is well defined over $[0,T]$ for all $R,T>0$. These are called {\em globally
convergent} series, and the set of all such series is denoted by
$\allseriesellGC$.
A Fliess operator $F_c$ defined on $B_{\mathfrak p}^m(R)[t_0,t_0+T]$
is said to be {\em realized}
by a state space realization
\begin{subequations}
\begin{align}
\dot{z}&= g_0(z)+\sum_{i=1}^m g_i(z)\,u_i,\;\;z(t_0)=z_0  \label{eq:state}\\
y&=h(z), \label{eq:output}
\end{align}
\end{subequations}
where each $g_i$ is an analytic vector field expressed in local
coordinates on some neighborhood ${\cal W}\subseteq \re^n$ of $z_0$,
and $h$ is an analytic function on ${\cal W}$, if
(\ref{eq:state}) has a well defined solution $z(t)$, $t\in[t_0,t_0+T]$ on ${\cal W}$
for any given input $u\in B_{\mathfrak
p}^m(R)[t_0,t_0+T]$, and
\begdi
F_c[u](t)=h(z(t)),\;\; t\in[t_0,t_0+T]. \label{eq:realization}
\enddi
In this case, the coefficients of the $i$-th component of generating series $c$ are computed by
\begeq
(c_i,\eta)=L_{g_{\eta}}h_i(z_0),\;\;\eta\in X^\ast, \label{eq:c-equals-Lgh}
\endeq
where
\begdi
L_{g_{\eta}}h_i:=L_{g_{j_1}}\cdots L_{g_{j_k}}h_i,
\;\;\eta=x_{j_k}\cdots x_{j_1},
\enddi
the {\em Lie derivative} of $h_i$ with respect to $g_j$ is defined as
\begdi
L_{g_j} h_i: W\rightarrow \re: z \mapsto \frac{\partial h_i}{\partial z}(z)\,
g_j(z),
\enddi
and $L_{g_\emptyset}h_i=h_i$.

When Fliess operators $F_c$ and $F_d$ are connected in a parallel-product
fashion, it is known that $F_cF_d=F_{c\shuffle d}$. If $F_c$ and $F_d$ with
$c\in\allseriesell$ and $d\in\allseriesm$ are interconnected in a cascade
manner, the composite system $F_c\circ F_d$ has the Fliess operator
representation $F_{c\circ d}$, where $c\circ d$ denotes the {\em composition
product} of $c$ and $d$ as described in \cite{Gray-etal_SCL14}.
This product is associative and $\re$-linear in its left argument $c$.
In the event that two Fliess operators are interconnected
to form a feedback system, the closed-loop system has a Fliess operator
representation whose generating series is the {\em feedback product} of $c$ and
$d$, denoted by $c@d$. This product can be explicitly computed via Hopf algebra
methods. The basic idea is to consider the set of operators
 $\Fliessdelta=\{I+F_c\;:\;c\in\allseriesm\}$,
where $I$ denotes the identity map, as a group under composition.
It is convenient to introduce the symbol $\delta$ as the (fictitious)
generating series for the identity map. That is, $F_\delta:=I$ such that
$I+F_c:=F_{\delta+c}=F_{c_\delta}$ with $c_\delta:=\delta+c$. The set of all
such generating series for $\Fliessdelta$ is denoted by $\allseriesdelta$.
This set also forms a group under the composition product induced by operator
composition, namely,  $c_\delta\circ d_\delta:=\delta+d+c\modcomp d$, where
$\modcomp$ denotes the {\em modified composition product}
\cite{Gray-etal_SCL14}. The group $(\allseriesdelta,\circ,\delta)$
has coordinate functions that form a Fa\`{a} di Bruno type Hopf algebra. In
which case, the group (composition) inverse $c_\delta^{\circ -1}$
can be computed efficiently via the antipode of this Hopf
algebra \cite{Duffaut-etal_JACxx,Gray-etal_MTNS14,Gray-etal_SCL14}. This
inverse also provide an explicit expression for 
the feedback product, namely, $c@d=c\circ(\delta-d\circ c)^{\circ -1}$.

\subsection{Left Inversion of Multivariable Fliess Operators}

It was shown in \cite{Wang_90} that $F_c$ will map every input which is
analytic at $t_0$ to an output which is also analytic at $t_0$ provided
$c\in\allseriesellLC$. In \cite{Gray-etal_CDC15} an explicit formula was
developed for calculating the left inverse of a multivariable mapping $F_c$ given a real
analytic function in its range. Without loss of generality assume $t_0=0$. Note
that every $c\in\allseries$ can be decomposed into its natural and forced
components, that is, $c=c_N+c_F$, where $c_N:=\sum_{k\geq 0} (c,x_0^k)x_0^k$
and $c_F:=c-c_N$. A condition under which the left inverse of $F_c$ exists is
provided by the following definition

\begde \label{de:vector-relative-degree-c}
Given $c\in\allseriesm$,
let $r_i\geq 1$ be the largest integer such that $\supp(c_{F,i})\subseteq
x_0^{r_i-1}X^\ast$,
where $i=1,2,\ldots,m$.
Then the component series $c_i$ has \bfem{relative degree} $r_i$ if the linear
word
$x_0^{r_i-1}x_j\in \supp(c_i)$ for
some $j\in\{1,\ldots,m\}$,
otherwise it is not well defined.
In addition, $c$ has \bfem{vector relative degree} $r=[r_1\;r_2\;\cdots\;r_m]$
if each $c_i$ has relative degree $r_i$ and the $m\times m$ matrix
\begdi
\footnotesize
A=\left[
\begin{tabular}{@{ }c@{$\;\;$ }c@{$\;\;$}c@{$\;\;$}c@{ }}
$(c_1,x_0^{r_1-1}x_1)$    &  $(c_1,x_0^{r_1-1}x_2)$    & $\cdots$  &
$(c_1,x_0^{r_1-1}x_m)$ \\
$(c_2,x_0^{r_2-1}x_1)$    &  $(c_2,x_0^{r_2-1}x_2)$    & $\cdots$  &
$(c_2,x_0^{r_2-1}x_m)$ \\
$\vdots$                  &  $\vdots$                  & $\vdots$  &  $\vdots$
\\
$(c_m,x_0^{r_m-1}x_1)$    &  $(c_m,x_0^{r_m-1}x_2)$    & $\cdots$  &
$(c_m,x_0^{r_m-1}x_m)$
\end{tabular}
\right]
\enddi
has full rank. Otherwise, $c$ does not have vector relative degree.
\endde

This notion of vector relative degree agrees with the usual definition
in the state space setting \cite{Isidori_95}. In particular, $c$ has vector
relative degree $r$ only if for each $i$ the series $(x_0^{r_i-1}x_j)^{-1}(c_i)$
is non proper for some $j$. Here the left-shift operator for any $x_i\in X$ is
defined on $X^\ast$ by $x_i^{-1}(x_i\eta)=\eta$ with $\eta\in X^\ast$ and zero
otherwise. Higher order shifts are defined inductively via
$(x_i\xi)^{-1}(\cdot)=\xi^{-1}x_i^{-1}(\cdot)$, where $\xi\in X^\ast$.
The left-shift operator is assumed to act linearly and componentwise on
$\allseriesm$. The shuffle inverse of any series $C\in\allseriesmtimesm$ is given
by
\begdi
C^{\shuffle -1}=((C,\emptyset)(I-C^\prime))^{\shuffle
-1}=(C^\prime)^{\shuffle \ast}(C,\emptyset)^{-1},
\enddi
where $C^\prime=I-(C,\emptyset)^{-1}C$ is proper, i.e.,
$(C^\prime,\emptyset)=0$, and $(C^\prime)^{\shuffle\ast}:=\sum_{k\geq 0}
(C^\prime)^{\shuffle k}$. The relationship between $C^{\shuffle -1}$ and the
multiplicative inverse operator $(F_C)^{-1}$, that is, $F_C(F_C)^{-1}=(F_C)^{-1}F_C=I$, is $(F_C)^{-1}=F_{C^{\shuffle -1}}$.

Let $X_0:=\{x_0\}$, and $\allseriesXO$ denotes the set of all commutative series
over $X_0$. When $c\in\allseriesXO$, $F_c[u](t)$ reduces to the Taylor series
$\sum_{k\geq 0} (c,x_0^k)E_{x_0^k}[u](t)=\sum_{k\geq 0} (c,x_0^k) t^k/k!$. The
main inversion tool used in the paper is given next.

\begth \label{th:main-left-inversion-formula}
Suppose $c\in\allseriesm$ has vector relative degree $r$. Let $y$ be analytic at
$t=0$
with generating series $c_y\in\allseriesXOmLC$ satisfying
$(c_{y_i},x_0^k)=(c_i,x_0^k)$, $k=0,1,\ldots,r_i-1$, $i=1,2,\ldots,m$.
Then the input
\begeq \label{u-Fdb-style1}
u(t)=\sum_{k=0}^\infty (c_u,x_0^k)\frac{t^k}{k!},
\endeq
is the unique real analytic solution to $F_c[u]=y$ on $[0,T]$ for some $T>0$,
where
\begeq \label{u-Fdb-style2}
c_u= \left(\left[C^{\shuffle -1}\shuffle (x_0^r)^{-1}(c-c_y)
\right]^{\circ -1}\right)_N,
\endeq
the $i$-th row of $(x_0^r)^{-1}(c-c_y)$ is $(x_0^{r_i})^{-1}(c_i-c_{y_i})$, and
the $(i,j)$-th entry of $C$ is $(x_0^{r_i-1}x_j)^{-1}(c_i)$.
\endth

Observe that this theorem describes the outputs that can be successfully tracked
by $F_c$, namely, those whose Taylor series coefficients satisfy certain matching
conditions. (The same conditions appear in the classical papers on dynamic inversion.)
For example, if $m=2$ and the vector relative degree is
$(2,2)$, the Taylor series of the output are subject to the
constraints $(c_{y_1},\emptyset)=(c_1,\emptyset)$,
$(c_{y_2},\emptyset)=(c_2,\emptyset)$, $(c_{y_1},x_0)=(c_1,x_0)$, and
$(c_{y_2},x_0)=(c_2,x_0)$, which in turn puts constraints on the class of 
admissible output paths. 

\section{Left Inverse of Bi-steerable Car Kinematics}
\label{sec:left-inversion-Bisteered-system}

\begin{figure}[t]
\begce
\includegraphics[width=8.5cm]{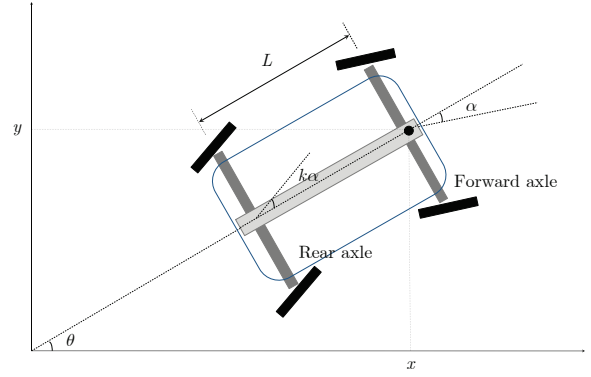}
\endce
\caption{Bi-steerable car}
\label{fig:bisteered_car}
\end{figure}

Consider the bi-steerable car shown in Figure~\ref{fig:bisteered_car}.
For simplicity, only the kinematics are considered. So the car is assumed to have
zero mass and move in
the plane with the speed of the car $u_1=\sqrt{\dot{x}^2+\dot{y}^2}$ and the front axle
steering angular velocity $u_2=\alpha$ as inputs.
The kinematics of the system are therefore
\begin{align*}
\dot{x} & =  u_1\cos (\theta+\alpha),\\
\dot{y} & =  u_1\sin (\theta+\alpha),\\
\dot{\theta} & =  u_1\frac{\sin(\alpha-f(\alpha))}{L\cos(f(\alpha))},\\
\dot{\alpha} & =  u_2,
\end{align*}
where $f(\alpha):=k\alpha$ with $k\in \re$. These dynamics assume the usual constraints of
rolling without slippage of the wheels, namely,
\begin{align*}
0 = & \dot{y}_{forward} \cos(\theta+\alpha)-\dot{x}_{forward} \sin(\theta +
\alpha),  \\
0 = & \dot{y}_{rear} \cos(\theta+f(\alpha))-\dot{x}_{rear} \sin(\theta +
f(\alpha)).
\end{align*}
Setting $z_1=x$, $z_2=y$, $z_3=\theta$, $z_4=\alpha$,
the outputs are picked to be the coordinates of
the front axle center $y_i= z_i$, $i=1,2$. The corresponding
two-input, two-output state space realization is
\begin{subequations} \label{eq:simplified_diff-system}
\begin{align}
\left(\begin{array}{c} \dot{z_1} \\ \dot{z_2} \\ \dot{z_3} \\
\dot{z_4}\end{array}\right) &
= \left(\begin{array}{c} \cos (z_3+z_4) \\ \sin(z_4+z_4) \\
\frac{\sin((1-k)\alpha)}{L \cos(k\alpha)} \\ 0 \end{array}\right) u_1 +
\left(\begin{array}{c} 0 \\ 0 \\ 0 \\ 1 \end{array}\right) u_2 \\
\left(\begin{array}{c} y_1 \\ y_2 \end{array}\right) & =
\left(\begin{array}{c} z_1 \\ z_2 \end{array}\right).
\end{align}
\end{subequations}
Hereafter, the focus is on the $k=-0.7$ case, which as explained in
\cite{Sekhavat-Hermosillo_02} is sufficient for the existence of flat outputs.
But this fact in inconsequential here since flat outputs will not be used.
The generating series,
$c$, can be computed directly from
\rref{eq:c-equals-Lgh}
using the vector fields and output function given in
\rref{eq:simplified_diff-system}:
\begin{align*}
\lefteqn{c_1 =  z_{1,0} + \co\, x_1 -\frac{\sin(1.7z_{4,0})}{\cos(0.7z_{4,0})}
\,\si\, x_1^2 -\si \,x_1x_2} \\
& -{\co}\, x_1x_2 x_2 -
\frac{\sin(1.7z_{4,0})}{\cos(0.7z_{4,0})} \,\co\, x_1x_2x_1 +\cdots\\
\lefteqn{c_2  =  z_{2,0} + \si\, x_1 +\frac{\sin(1.7z_{4,0})}{\cos(0.7z_{4,0})}
\,\co\, x_1^2 +\co\, x_1x_2} \\
& -{\si}\, x_1x_2 x_2 -
\frac{\sin(1.7z_{4,0})}{\cos(0.7z_{4,0})} \,\si\, x_1x_2x_1
+\cdots,
\end{align*}
where $\co:=\cos(z_{3,0}+z_{4,0})$ and $\si:=\sin(z_{3,0}+z_{4,0})$.
\begin{figure}[t]
\begce
\includegraphics[width=9.cm]{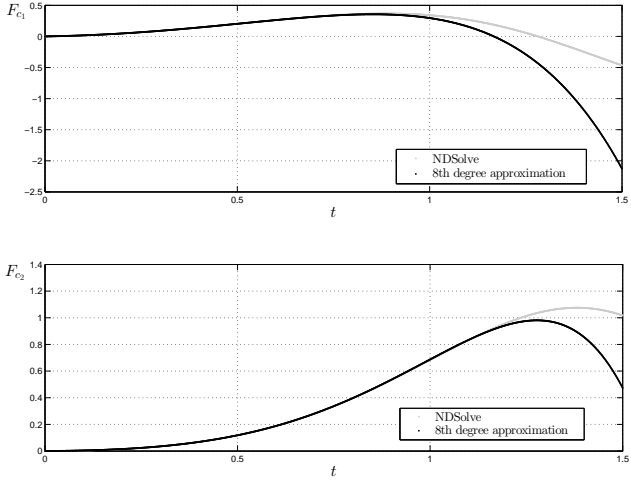}
\endce
\caption{Mathematica computed (NDSolve) and Fliess operator approximated (degree 8) step response of bi-steerable car ($u_1(t)=1.5$ and
$u_2(t)=0.8$, $t>0$)}
\label{fig:bisteered_car_sim}
\end{figure}%
In this case, $c\in\re_{GC}^2\langle
\langle X \rangle\rangle$ with $X=\{x_0,x_1,x_2\}$ since all sines and cosines 
in the numerator of every coefficient can be bounded by $1$. Therefore, after 
some additional calculations a conservative choice for the growth constants is 
$K_c= \max\{z_{1,0},z_{2,0}\}$ and $M_c=2.4\sec(0.7 z_{4,0})$. 
With $z_{1,0} =
z_{2,0}=z_{3,0}=0$, $z_{4,0}=z_{5,0}=0.2$, the step response ($u_1(t)=1.5$ and $u_2(t)=0.8$, $t>0$)
of the system was computed using the {\tt NDSolve} Mathematica function
\cite{Mathematica} and is compared against a degree eight Fliess operator
approximation in
Figure~\ref{fig:bisteered_car_sim}.  Since $c\in \re_{GC}^2\langle
\langle X \rangle\rangle$, the output $y=F_c[u]$ is well defined on
$[0,\infty)$ for all $u\in B_1^m[0,T](R)$ and any finite $R,T>0$. \cite{Winter-Gray-Duffaut_Espinosa_MTNS2016}.

A simple calculation shows that the relative degrees of $y_1$ and $y_2$ are both
one, but the system does not have a well defined vector relative degree.
Applying dynamic extension on $u_1$ yields the
augmented system
\begin{align} \label{eq:augmented_car_model}
 \left(\begin{array}{c} \dot{z_1} \\ \dot{z_2} \\ \dot{z_3} \\
\dot{z_4}\\ \dot{z_5}\end{array}\right)
\!=\! \left(\begin{array}{c} z_5\cos (z_3+z_4) \\ z_5\sin(z_4+z_4) \\
z_5\frac{\sin((1-k)\alpha)}{L \cos(k\alpha)} \\ 0 \\ 0 \end{array}\right)
\!+\! \left(\begin{array}{c} 0 \\ 0 \\ 0 \\ 1 \\ 0 \end{array}\right)u_2 \!+\!
\left(\begin{array}{c} 0 \\ 0 \\ 0 \\ 0 \\ 1 \end{array}\right) \bar{u}_1,
\end{align}
where the new input $\bar{u}_1$ is taken as the derivative of the original
input $u_1$. This system has
vector relative degree $r=[r_1\;r_2]=[2\;2]$ with $r_1+r_2=4<5=n$ and
generating series
\begin{align*}
c_1 = & \, z_{1,0} + z_{5,0}\,\co \, x_0 + z_{5,0}\,\si \,x_0x_1  + \co\, x_0
x_2\\
&- z_{5,0}^2\frac{\sin(1.7z_{4,0})}{\cos(0.7z_{4,0})}\,\si\, x_0^2 +\cdots\\
c_2  =& \, z_{2,0} + z_{5,0}\, \si\,  x_0 + z_{5,0}\, \co\,  x_0x_1 + \si\, x_0
x_2\\
& + z_{5,0}^2\frac{\sin(1.7z_{4,0})}{\cos(0.7z_{4,0})}\,\co \, x_0^2 +\cdots.
\end{align*}
These series are also elements of $\re_{GC}^2\langle\langle X \rangle\rangle$
since again all sines and cosines in the numerator can be 
bounded by $1$ so that conservative growth constants are
$K_c=\max\{z_{1,0},z_{2,0},z_{5,0}\}$ and $M_c=2.4 z_{5,0}\sec(0.7 
z_{4,0})$. The decoupling matrix
\begeq \label{eq:dec_matr_A}
A=(C,\emptyset) = \left( \begin{array}{cc} z_{5,0}\,\si &
\co \\ z_{5,0}\,\co & \si
\end{array}
  \right)
\endeq
is clearly nonsingular as long as $z_{5,0}\neq 0$.

Given a desired output function
\begin{align} \label{eq:cyFPS}
y(t) & = \sum_{k=0}^\infty (c_y,x_0^k)\frac{t^k}{k!},
\end{align}
where $c_y = [c_{y_1},c_{y_2}]^T$ is the generating series of
$y$, the left inverse $c_u=[c_{\bar{u}_1},c_{u_2}]^T $ is computed directly from
\rref{u-Fdb-style1}-\rref{u-Fdb-style2}.  It is sufficient here to consider polynomial outputs up to degree three, so
let $(c_{y_j},x_0^i) = v_{ij}$ for $i=0, 1,2,3$
and $j=1,2$. The series
\begdi d=\left(\begin{array}{c} d_1\\
d_2  \end{array} \right) :=C^{\shuffle -1}\shuffle (x_0^r)^{-1}(c-c_y)
\enddi
is found to be
\begin{align*}
\lefteqn{d_1 = \frac{v_{12}\,\si}{z_{5,0}}
-\frac{v_{22}\,\co}{z_{5,0}} + z_{5,0} \frac{\sin(1.7
z_{4,0})}{\cos(0.7z_{4,0})} } \\
& + \left( \left(\frac{ v_{13}  }{z_{5,0}} +v_{22} \frac{\sin(1.7
z_{4,0})}{\cos(0.7z_{4,0})}  \right)\,\si \right. \\
& \left.+\left(v_{12}\frac{\sin(1.7 z_{4,0})}{\cos(0.7z_{4,0})}-
\frac{v_{23}}{z_{5,0}}  \right)\,\co \right)x_0 \\
& + \left(\frac{v_{22}\,\si}{z_{5,0}}  +1.7 z_{5,0} \frac{\cos(1.7
z_{4,0})}{\cos(0.7z_{4,0})} \right.\\
& \left. +  0.7 z_{5,0}\frac{\sin(1.7
z_{4,0})\tan(0.7z_{4,0})}{\cos(0.7z_{4,0})}
+\frac{v_{22}\,\co}{z_{5,0}} \right)x_1 \\
& + \left(\frac{v_{22}\,\co}{z_{5,0}^2}
-\frac{v_{12}\,\co}{z_{5,0}^2} + \frac{\sin(1.7 z_{4,0})}{\cos(0.7z_{4,0})}
\right)x_2 +\cdots,\\
\lefteqn{d_2  = -v_{12}\,\co -v_{22}\,\si } \\
& - \left(v_{13}\co -z_{5,0} \frac{\sin(1.7 z_{4,0})}{\cos(0.7z_{4,0})}
\left(v_{22}\,\co+v_{12}\, \si \right)
+ v_{23}\si \right)x_0 \\
& + \left(v_{12}\,\si - v_{22}\,\co \right)x_1 +\cdots
\end{align*}
The composition inverse of $d$ is computed componentwise using the recursive
method in \cite{Duffaut-etal_JACxx,Gray-etal_MTNS14,Gray-etal_SCL14}.  In which
case, the formulas for $c_{\bar{u}_1}$ and $c_{u_2}$ are, respectively,
\begin{align*}
\lefteqn{c_{\bar{u}_1}=\left(d_1^{\circ{-1}}\right)_N= \frac{v_{22} 
\si}{z_{5,0}}
-\frac{v_{12} \co}{z_{5,0}} - z_{5,0}\frac{\sin(1.7
z_{4,0})}{\cos(0.7z_{4,0})}}\\
& + \frac{1}{z_{5,0}^2} \left(\rule{0.in}{0.2in} -2 v_{12} v_{22}
\cos(2(z_{3,0}+z_{4,0})) \right. \\
& +1.7 z_{5,0}^4 \frac{\sin(1.7z_{4,0})\cos(1.7
z_{4,0})}{\cos(0.7z_{4,0})^2} - v_{13} z_{5,0} \,\si \\
& + 1.7 v_{12} z_{5,0}^2
\frac{\cos(1.7 z_{4,0})}{\cos(0.7z_{4,0})} \,\si - v_{22} z_{5,0}^2
\frac{\sin(1.7 z_{4,0})}{\cos(0.7z_{4,0})}\, \si \\
& + v_{12}^2 \sin(2(z_{3,0}+z_{4,0})) - v_{22}^2
\sin(2(z_{3,0}+z_{4,0}))\\
& +0.7 z_{5,0}^4 \frac{\sin(1.7 z_{4,0})^2\tan(0.7
z_{4,0})}{\cos(0.7z_{4,0})^2} \\
&  + z_{0,5} \,\co \, \left( v_{23} - 1.7 v_{22}
z_{5,0} \frac{\cos(1.7 z_{4,0})}{\cos(0.7z_{4,0})} \right. \\
& \left.\left. - z_{5,0}\frac{\sin(1.7
z_{4,0})}{\cos(0.7z_{4,0})} (v_{12}+0.7 v_{22} \tan(0.7 z_{4,0})) \right)
\right)x_0 +\cdots\\
\lefteqn{c_{u_2}=\left(d_2^{\circ{-1}}\right)_N= v_{22} \si + v_{12} \co  +
\frac{1}{z_{5,0}} \left( \frac{v_{12}^2}{2} + \frac{v_{22}^2}{2} \right. } \\
& + v_{13} z_{5,0} \, \co +
\frac{(v_{22}^2-v_{12}^2)}{2}\cos(2(z_{3,0}+z_{4,0})) + v_{23}z_{5,0} \, \si \\
& \left. - v_{12} v_{22} \sin(2(z_{3,0}+z_{4,0})) \rule{0.in}{0.2in}\right) x_0
+\cdots
\end{align*}
Some key features of these expressions are:
\begin{enumerate}

\item[$i.$] The formulas are exact if not truncated. But, of course,
truncation is necessary for implementation.
The truncation gives the degree of approximation. Here the degree of 
approximation for $c_{\bar{u}_1}=(d_1^{\circ-1})_N$ and 
$c_{u_2}=(d_2^{\circ-1})_N$ means their truncation is to degree $6$. These 
truncated inputs are then fed into the $8$th degree approximation of series $c$ 
via the composition product described in \cite{Gray-etal_SCL14}.


 \item[$ii.$] The formulas only need to be computed once. They can therefore
be loaded directly into a micro-controller inside the vehicle. Then
based on measurements of the current position and steering angle, the input for
tracking the next section of the desired path can be quickly computed by just
numerically evaluating the formula.

\item[$iii.$] Flat outputs are not required for
computing these left inverses.

\item[$iv.$] One can increase the degree of approximation of the output tracking
by including more input terms if
the computing power is available.

\end{enumerate}

\section{Collision Avoidance System} \label{sec:simulations}

\begin{figure}[t]
\begin{center}
\includegraphics[width=8.cm]{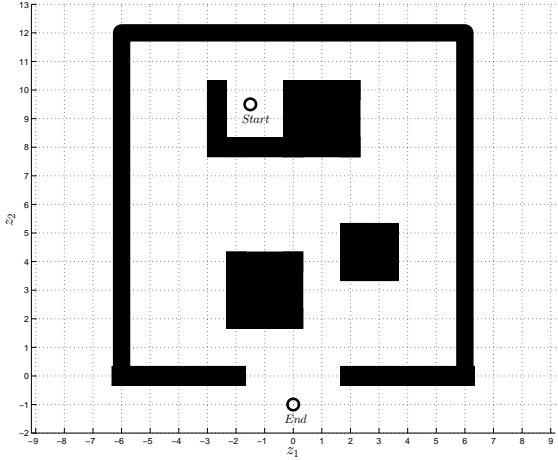}
\vspace*{-0.1in}
\caption{Map of area with obstacles}
\vspace*{-0.1in}
\label{fig:area}
\end{center}
\end{figure}

\begin{figure}[t]
\begin{center}
\includegraphics[width=8.cm]{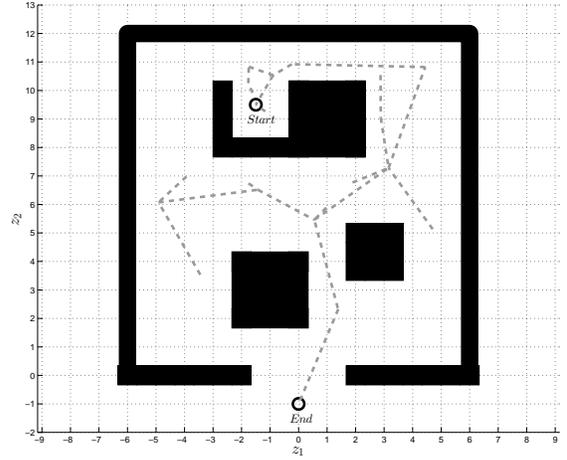}
\vspace*{-0.1in}
\caption{RRT used to identify the shortest path}
\vspace*{-0.1in}
\label{fig:tree}
\end{center}
\end{figure}

\begin{figure}[t]
\begin{center}
\includegraphics[width=8.cm]{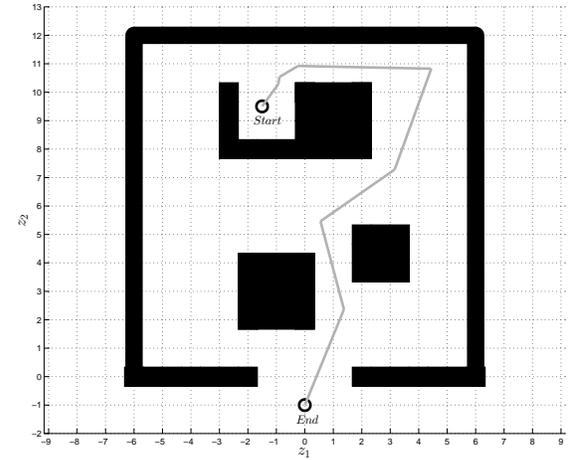}
\vspace*{-0.1in}
\caption{Shorted path extracted from the RRT}
\vspace*{-0.1in}
\label{fig:tree_path}
\end{center}
\end{figure}


A collision avoidance system is now described based on the
left inversion input formula for the bi-steerable car model developed in the previous
section. The basic steps of the algorithm are as follows:

\begin{enumerate}
\item[$i.$] Map the area and obstacles where the bi-steerable car moves
and provide a \emph{start} location and an \emph{end} location.

\item[$ii.$] 
The RRT path planning algorithm described in
\cite{Lavalle_2006} then uses the data from the previous step
to compute the path between the start and end locations by first computing a
tree and then identifying the shortest path within the tree
from the start location to the end location.

\item[$iii.$] The shortest path is then smoothed and approximated by a cubic spline.

\item[$iv.$] Each spline section, which is a cubic polynomial, is fed to the left inversion formula in order
to compute the corresponding control inputs $\bar{u}_1$ and $u_2$.

\item[$v.$] The inputs $\bar{u}_1$ and $u_2$ are finally used to drive the bi-steerable car.

\end{enumerate}

\vspace*{0.1in}

Consider the following specific example: \\

\noindent{\bf\em i}: Suppose the area with obstacles is the one shown in
Figure~\ref{fig:area}. Pick the start and end locations to be:
\begin{align*}
Start:& \;\;\;(-1.5,9.5) \\
End: & \;\;\;(0,-1).
\end{align*} 

\noindent{\bf\em ii}: The RRT generated for this map
is shown in
Figure~\ref{fig:tree}. The shortest path between the 
given start and end locations exacted from this RRT is shown 
in Figure~\ref{fig:tree_path}. \\

\begin{figure}[t]
\begin{center}
\includegraphics[width=8.cm]{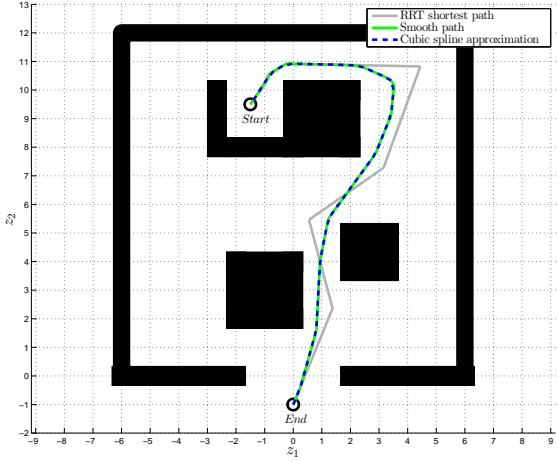}
\vspace*{-0.1in}
\caption{Shortest, smoothed and spline approximated paths}
\vspace*{-0.1in}
\label{fig:smoothed_path}
\end{center}
\end{figure}

\noindent{\bf\em iii}: The shortest path is then smoothed and approximated by cubic
splines as shown in Figure~\ref{fig:smoothed_path}. The constant terms of each
cubic spline section must coincide with the final position of the previous section,
This means that one does not need
to fit the constant terms of the spline approximation. That is, the
condition $(c,\emptyset)=(c_y,\emptyset)$ in Theorem~\ref{th:main-left-inversion-formula} 
is always satisfied by the cubic splines in
the path approximations. The linear coefficient of the spline approximation is also
subject to the range conditions of
Theorem~\ref{th:main-left-inversion-formula}. However,
since $z_4$ and $z_5$ in \rref{eq:augmented_car_model} are
just the integrals of $u_2$ and $\bar{u}_1$, respectively, their
corresponding initial conditions $z_{4,0}$ and $z_{5,0}$ can be chosen
arbitrarily as long as $z_{5,0}\neq 0$, otherwise the decoupling matrix
\rref{eq:dec_matr_A} is singular. Observe that
\begeq \label{eq:first_order_coeffs}
(c,x_0) = \left( \begin{array}{c} z_{5,0} \cos(z_{3,0}+z_{4,0}) \\ z_{5,0}
\sin(z_{3,0}+z_{4,0})   \end{array} \right),
\endeq
where $z_{3,0}$ is known and fixed. Hence, one can always find $z_{4,0}$ and
$z_{5,0}$ such that $(c,x_0)$ matches the desired coefficient $(c_y,x_0)$,
which comes from the spline approximation. Here the time interval for the
simulation have been normalized to $[0,1]$ since only the kinematics are 
being considered.\\

\begin{figure}[t]
\begin{center}
\includegraphics[width=9.0cm]{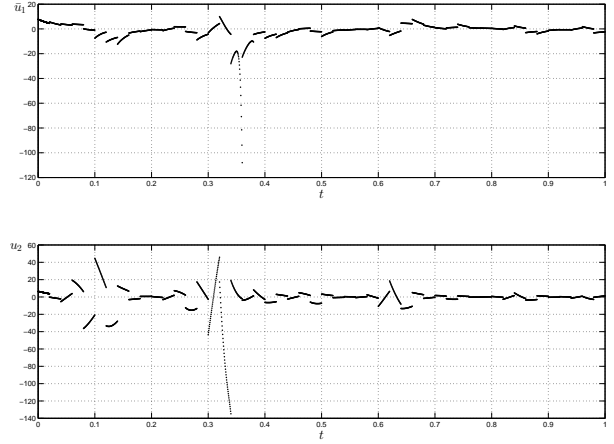}
\vspace*{-0.1in}
\caption{Control inputs $\bar{u}_{1}$ and $u_2$} 
\vspace*{-0.1in}
\label{fig:u1u2}
\end{center}
\end{figure}

\noindent{\bf\em iv}: The control inputs are now computed section by section. 
The first
spline section has initial conditions
\begdi
z_{1,0}=-1.5, \;\;z_{2,0} = 9.5, \;\; z_{3,0} = 1.5,
\enddi
and the desired outputs of the path are
\begdi
c_{y_1} = -1.5 + 10.75 x_0 - 28.15 x_0^2 + 185.15 x_0^3
\enddi
and
\begdi
c_{y_2} = 9.5 + 14.40 x_0 + 13.07 x_0^2 - 24.89 x_0^3.
\enddi
As expected, the constant coefficients are automatically equal to the start
position of the car. Also, \rref{eq:first_order_coeffs} is solved for $z_{4,0}$
and $z_{5,0}\neq 0$ so that $(c_1,x_0) = (c_{y_1}, x_0)=10.75$ and
$(c_2,x_0) = (c_{y_2}, x_0)=14.40$, which gives
\begdi
z_{4,0}= -4.00028,\;\; z_{5,0} = -17.97.
\enddi
This solution is not unique. Since the overall time of the simulation was
normalized, the time interval for each section of the path was chosen to have
duration $t_{sim}= 0.02$. The left inverse formulas in this case give
\begin{align*}
\bar{u}_1(t)= & 7.74 - 195.52 t + 1592.16 t^2 + 12925.9 t^3 \\
& - 443203.0 t^4 +  4.26\times 10^6 t^5 + 2.83\times 10^6 t^6
\end{align*}
and
\begin{align*}
u_2(t) = &\, 6.36 - 142.17 t + 386.62 t^2 - 226.02 t^3 \\
& - 2648.52 t^4 + 1280.47 t^5 + 76475.8 t^6.
\end{align*}
Since the degree of the approximation was chosen to be $6$, these inputs give 
the following errors:
\begin{align*}
y_1(t)-\hat{y}_1(t) = & {\;} 1.33\times 10^6 t^7 - 5.10\times 10^6 t^8 \\
& - 2.02\times 10^7 t^9 +  3.47\times  10^7 t^{10}+\cdots
\end{align*}
and
\begin{align*}
y_2(t)-\hat{y}_2(t) = & {\;} -1.24\times 10^6 t^7 + 4.39\times 10^6 t^8 \\
&+  2.09\times 10^7 t^9 -  2.82\times 10^7 t^{10}+\cdots,
\end{align*}
where $y_i$ is the $i$-th component of $y$ having generating series $c_y$, and
$\hat{y}_i = c\circ c_{u}$ is truncated to degree $6$.
The computed control inputs are shown in Figure~\ref{fig:u1u2}.

\begin{figure}[h]
\begin{center}
\includegraphics[width=8.cm]{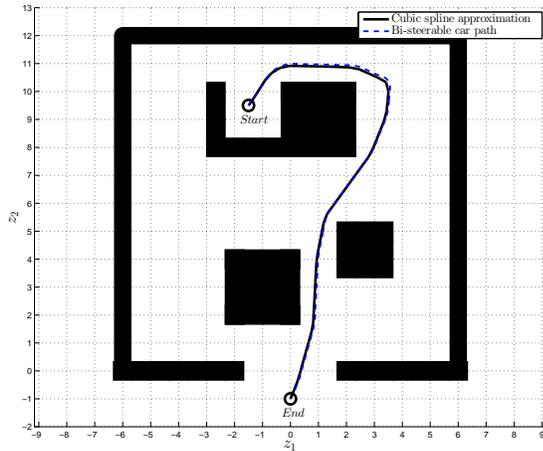}
\vspace*{-0.1in}
\caption{Cubic spline approximated path and bi-steerable car trajectory} 
%
\vspace*{-0.1in}
\label{fig:generated_path}
\end{center}
\end{figure}

\noindent{\bf\em iv}: Finally, the computed control inputs
$\bar{u}_1$ and $u_2$ are used to drive the bi-steerable car as shown in
Figure~\ref{fig:generated_path}. In the same figure, the 
bi-steerable car path is overlapped, as comparison, with the cubic spline 
approximation of the shortest path (avoiding obstacles) computed by the RRT 
algorithm. The degree of approximation and smoothness of the bi-steerable car 
path can be tuned by increasing or decreasing the number of partitions of the 
cubic spline approximation, and the degree of approximation of the computed left 
inverses. Tracking performance will ultimately be bounded by the amount of 
computational power available. 

\section{Conclusions} \label{sec:comclusions}

A collision avoidance system was described for a bi-steerable car
based on a left inversion formula for Fliess operators whose generating
series have a well defined vector relative degree. This allows one to
integrate directly the kinematics of the car and a path planning algorithm
without the need for passing through a flat output as is done in other
methods.
In principle, the full dynamics could also be inverted to give a more
realistic collision avoidance system.
The method was demonstrated by numerical simulation.

\end{document}